\documentclass[12pt,a4paper]{article}
\usepackage{amsmath}
\usepackage{amssymb}
\usepackage{enumerate}
\usepackage{t1enc}
\usepackage[latin1]{inputenc}
\usepackage[english]{babel}
\usepackage{enumerate}

\usepackage{mathtools}

\usepackage{amsfonts}
\usepackage{latexsym}
\usepackage{bm}
\usepackage{tikz}

\newtheorem{theorem}{Theorem} 
\newtheorem{prop}{Proposition}
\newtheorem{lemma}{Lemma}
\newtheorem{cor}{Corollary}

\newtheorem{construction}{Construction}

\allowdisplaybreaks

\begin{document}
\title{Isolation of non-triangle cycles in graphs}

\author{Peter Borg\footnote{Email address: peter.borg@um.edu.mt} \quad \quad \quad Dayle Scicluna\footnote{Email address: dayle.scicluna.09@um.edu.mt}\\[5mm]
{\normalsize Department of Mathematics} \\
{\normalsize Faculty of Science} \\
{\normalsize University of Malta}\\
{\normalsize Malta}\\
}

\date{}
\maketitle

\begin{abstract} 
Given a set $\mathcal{F}$ of graphs, we call a copy of a graph in $\mathcal{F}$ an $\mathcal{F}$-graph. The $\mathcal{F}$-isolation number of a graph $G$, denoted by $\iota(G, \mathcal{F})$, is the size of a smallest set $D$ of vertices of $G$ such that the closed neighbourhood $N[D]$ of $D$ intersects the vertex sets of the $\mathcal{F}$-graphs contained by $G$ (equivalently, $G-N[D]$ contains no $\mathcal{F}$-graph). Let $\mathcal{C}$ be the set of cycles, and let $\mathcal{C}'$ be the set of non-triangle cycles (that is, cycles of length at least $4$). Let $G$ be a connected graph having exactly $n$ vertices and $m$ edges. The first author proved that $\iota(G,\mathcal{C}) \leq n/4$ if $G$ is not a triangle. Bartolo and the authors determined a set $\mathcal{E}_4$ of nine graphs such that $\iota(G,\{C_4\}) \leq n/5$ if $G$ is not an $\mathcal{E}_4$-graph. Chen and Cui recently strengthened this result by proving that $\iota(G, \mathcal{C}') \leq n/5$ if $G$ is not an $\mathcal{E}_4$-graph. This had been proved by Zhang and Wu for the case where $G$ contains no triangle. Various authors proved that $\iota(G,\mathcal{C}) \leq (m+1)/5$ if $G$ is not a triangle. We prove that $\iota(G,\mathcal{C}') \leq (m+1)/6$ if $G$ is not a $4$-cycle. Zhang and Wu established this for the case where $G$ contains no triangle. Our result yields the inequality $\iota(G,\{C_4\}) \leq (m+1)/6$ of Wei, Zhang and Zhao. These bounds are attained by infinitely many (non-isomorphic) graphs. The proof of our inequality hinges on also determining the graphs attaining the bound.
\end{abstract}

\section{Introduction}
The study of isolating sets was initiated by Caro and Hansberg~\cite{CaHa17}. It generalizes the study of the classical domination problem \cite{C, CH, HHS, HHS2, HL, HL2} naturally. Indeed, $D$ is a \emph{dominating set of $G$} (that is, the closed neighbourhood of $D$ is the vertex set of $G$) if and only if $D$ is a $K_1$-isolating set of $G$, so the \emph{domination number} is the $K_1$-isolation number. Isolation of cycles has attracted particular attention. A cycle isolation problem posed in \cite{CaHa17} was solved by the present first author in \cite{Borg1}, and this motivated further investigations that yielded the results in \cite{BBS, CC, CCZ, CZ, WZZ, Y, ZW2, ZW3, ZW4} among others. In this paper, we solve a problem that arises from these works.

Unless stated otherwise, we use small letters such as $x$ to denote non-negative integers or elements of sets, and capital letters such as $X$ to denote sets or graphs. For $n \geq 0$, $[n]$ denotes the set $\{i \in \mathbb{N} \colon i \leq n\}$, where $\mathbb{N}$ is the set of positive integers. Note that $[0]$ is the empty set $\emptyset$. Arbitrary sets are taken to be finite. For a set $X$, ${X \choose 2}$ denotes the set of $2$-element subsets of $X$. We may represent a $2$-element set $\{x,y\}$ by $xy$.

For standard terminology in graph theory, we refer the reader to \cite{West}. Most of the notation and terminology used here is defined in \cite{Borg1}, which motivates the work in this paper. 

Every graph $G$ is taken to be \emph{simple}, that is, $G$ is a pair $(V(G), E(G))$ such that $V(G)$ and $E(G)$ (the vertex set and the edge set of $G$) are sets that satisfy $E(G) \subseteq {V(G) \choose 2}$. We call $G$ an \emph{$n$-vertex graph} if $|V(G)| = n$. We call $G$ an \emph{$m$-edge graph} if $|E(G)| = m$. For a vertex $v$ of $G$, $N_{G}(v)$ denotes the set of neighbours of $v$ in $G$, $N_{G}[v]$ denotes the closed neighbourhood $N_{G}(v) \cup \{v\}$ of $v$, and $d_{G}(v)$ denotes the degree $|N_{G} (v)|$ of $v$. For a subset $X$ of $V(G)$, $N_G[X]$ denotes the closed neighbourhood $\bigcup_{v \in X} N_G[v]$ of $X$, $G[X]$ denotes $(X, E(G) \cap {X \choose 2})$ (the subgraph of $G$ induced by $X$), and $G - X$ denotes the graph $G[V(G) \backslash X]$ (obtained by deleting the vertices in $X$ from $G$). We may abbreviate $G-\{x\}$ to $G-x$. Where no confusion arises, the subscript $G$ may be omitted from any notation that uses it; for example, $N_G(v)$ may be abbreviated to $N(v)$. If $H$ is a subgraph of $G$, then we say that \emph{$G$ contains $H$}. If $F$ is a copy of $G$, then we write $F \simeq G$. 

For $n \geq 1$, the graphs $([n], {[n] \choose 2})$ and $([n], \{\{i,i+1\} \colon i \in [n-1]\})$ are denoted by $K_n$ and $P_n$, respectively. For $n \geq 3$, $C_n$ denotes the graph $([n], \{\{1,2\}, \{2,3\}, \dots, \{n-1,n\}, \{n,1\}\})$. A copy of $K_n$ is called an \emph{$n$-clique} or a \emph{complete graph}. A copy of $P_n$ is called an \emph{$n$-path} or simply a \emph{path}. A copy of $C_n$ is called an \emph{$n$-cycle} or simply a \emph{cycle}. A $3$-clique is a $3$-cycle and is also called a \emph{triangle}. 

If $\mathcal{F}$ is a set of graphs and $F$ is a copy of a graph in $\mathcal{F}$, then we call $F$ an \emph{$\mathcal{F}$-graph}. A subset $D$ of $V(G)$ is called an \emph{$\mathcal{F}$-isolating set of $G$} if $N[D]$ intersects the vertex sets of the $\mathcal{F}$-graphs contained by $G$. Thus, $D$ is an $\mathcal{F}$-isolating set of $G$ if and only if $G-N[D]$ contains no $\mathcal{F}$-graph. It is to be assumed that $(\emptyset, \emptyset) \notin \mathcal{F}$. The size of a smallest $\mathcal{F}$-isolating set of $G$ is denoted by $\iota(G, \mathcal{F})$ and is called the \emph{$\mathcal{F}$-isolation number of $G$}. If $\mathcal{F} = \{F\}$, then we may replace $\mathcal{F}$ in these defined terms and notation by $F$.

One of the earliest domination results is the upper bound $n/2$ of Ore \cite{Ore} on the domination number of any connected $n$-vertex graph $G \not\simeq K_1$ (see \cite{HHS}). While deleting the closed neighbourhood of a dominating set produces the graph with no vertices, deleting the closed neighbourhood of a $K_2$-isolating set produces a graph with no edges. In the literature, a $K_2$-isolating set is also called a \emph{vertex-edge dominating set}. Caro and Hansberg~\cite{CaHa17} proved that if $G$ is a connected $n$-vertex graph, then $\iota(G, K_2) \leq n/3$ unless $G \simeq K_2$ or $G \simeq C_5$. This was independently proved by \.{Z}yli\'{n}ski \cite{Z} and solved a problem in \cite{BCHH}. Fenech, Kaemawichanurat and the present first author \cite{BFK} generalized these bounds by showing that for any $k \geq 1$, $\iota(G, K_k) \leq n/(k+1)$ unless $G \simeq K_k$ or $k=2$ and $G \simeq C_5$. This sharp bound settled a problem of Caro and Hansberg~\cite{CaHa17}. The graphs attaining the bound are determined in \cite{BG, CCZ, CCZ2, LMS}. Fenech, Kaemawichanurat and the present first author \cite{BFK2} also showed that $\iota(G, K_k) \leq (m+1)/({k \choose 2}+2)$ unless $G \simeq K_k$, and they determined the graphs attaining the bound. Generalizations of these bounds are given in \cite{Borg2, Borg4}.

Let $\mathcal{C}$ be the set of cycles. The present first author \cite{Borg1} obtained the following bound on $\iota(G, \mathcal{C})$, and consequently settled another problem of Caro and Hansberg~\cite{CaHa17}. 
\begin{theorem}[\cite{Borg1}] \label{Borgcycle}
If $G$ is a connected $n$-vertex graph that is not a triangle, then
\[\iota(G, \mathcal{C}) \leq \frac{n}{4}.\] 
Moreover, the bound is sharp.
\end{theorem}
He also gave an explicit construction of a connected $n$-vertex graph 
that attains the bound $\lfloor n/4 \rfloor$ resulting from Theorem~\ref{Borgcycle}. The graphs that attain the bound $n/4$ in the theorem are determined in \cite{CCZ}. Various authors have obtained an analogue of Theorem~\ref{Borgcycle} that provides a sharp bound on $\iota(G, \mathcal{C})$ in terms of the number of edges of $G$. In order to state the full result, we need the following construction \cite[Construction~1]{Borg3}, which is a generalization of \cite[Construction 1.2]{BFK2} and a slight variation of the construction in \cite[Section~1]{Borg1}.

\begin{construction}[\cite{Borg3}] \label{const} \emph{Consider any $m, k \in \{0\} \cup \mathbb{N}$ and any connected $k$-edge graph $F$, where $F \simeq K_1$ if $k = 0$ (that is, $V(F) \neq \emptyset$). By the division algorithm, there exist $q, r \in \{0\} \cup \mathbb{N}$ such that $m+1 = q(k+2) + r$ and $0 \leq r \leq k+1$. Let $Q_{m,k}$ be a set of size $q$. If $q \geq 1$, then let $v_1, \dots, v_q$ be the elements of $Q_{m,k}$, let $F_1, \dots, F_q$ be copies of $F$ such that the $q+1$ sets $V(F_1), \dots, V(F_q)$ and $Q_{m,k}$ are pairwise disjoint, and for each $i \in [q]$, let $w_i \in V(F_i)$, and let $G_i$ be the graph with $V(G_i) = \{v_i\} \cup V(F_i)$ and $E(G_i) = \{v_iw_i\} \cup E(F_i)$. If either $q = 0$, $T$ is the null graph $(\emptyset, \emptyset)$, and $G$ is a connected $m$-edge graph $T'$, or $q \geq 1$, $T$ is a tree with vertex set $Q_{m,k}$ (so $|E(T)| = q-1$), $T'$ is a connected $r$-edge graph with $V(T') \cap \bigcup_{i=1}^q V(G_i) = \{v_q\}$, and $G$ is a graph with $V(G) = V(T') \cup \bigcup_{i=1}^q V(G_i)$ and $E(G) = E(T) \cup E(T') \cup \bigcup_{i=1}^q E(G_i)$, then we say that $G$ is an \emph{$(m,F)$-special graph} with \emph{quotient graph $T$} and \emph{remainder graph $T'$}, and for each $i \in [q]$, we call $G_i$ an \emph{$F$-constituent of $G$}, and we call $v_i$ the \emph{$F$-connection of $G_i$ in $G$}. We say that an $(m,F)$-special graph is \emph{pure} if its remainder graph has no edges (see Figure~1). Clearly, an $(m,F)$-special graph is a connected $m$-edge graph.}
\end{construction}

\begin{center}
\begin{tikzpicture}[scale=1.05,
    vertex/.style={fill, circle, inner sep=0pt, minimum size=2.5pt}
]

    \foreach \i [evaluate=\i as \x using {(\i-1)*1.9}] in {1,...,7} {
        
        \node[vertex] (t\i) at (\x, 3.0) {};
        
        \node[vertex] (p\i) at (\x, 2.0) {};
        
        \node[vertex, label={below:$G_\i$}] (b\i) at (\x, 0.6) {};
        
        \node[vertex] (v\i L) at (\x-0.7, 1.3) {};
        \node[vertex] (v\i R) at (\x+0.7, 1.3) {};

        \draw[thick] (t\i.center) -- (p\i.center);
        
        \draw[thick] (p\i.center) -- (v\i L.center) -- (b\i.center) -- (v\i R.center) -- cycle;
    }

    \draw[thick, dotted] (t1) -- (t2);
    \draw[thick, dotted] (t2) -- (t3);
    \draw[thick, dotted] (t2) .. controls (3.6,3.25) and (4.0,3.25) .. (t4);
    \draw[thick, dotted] (t2) .. controls (4.2,3.6) and (5.2,3.6) .. (t5);
    \draw[thick, dotted] (t4) .. controls (7.2,3.25) and (7.6,3.25) .. (t6);
    
    \draw[thick, dotted] (t3) .. controls (5.5,3.7) and (9.5,3.7) .. (t7);

\end{tikzpicture}
\\
{\footnotesize \textbf{Figure 1:} An illustration of a pure $(m,C_4)$-special graph with $m = 41$.}
\end{center}

\begin{theorem}[\cite{Borg4,CZ,ZW2}] \label{Borgcycle2}
If $G$ is a connected $m$-edge graph that is not a triangle, then
\[\iota(G, \mathcal{C}) \leq \frac{m+1}{5}.\] 
Moreover, equality holds if and only if $G$ is a pure $(m,C_3)$-special graph or a $4$-cycle.
\end{theorem}

Theorem~\ref{Borgcycle} has inspired many other results. Consider a connected graph $G$, and let $n = |V(G)|$ and $m = |E(G)|$. Let $\mathcal{C}'$ be the set of cycles that are not triangles. Thus, $\mathcal{C}' = \{H \in \mathcal{C} \colon |V(H)| \geq 4\}$ and $C_4 \in \mathcal{C}'$. Strengthening a result in \cite{Y}, Bartolo and the present authors \cite{BBS} determined a set $\mathcal{E}_4$ of nine graphs such that $\iota(G,\{C_4\}) \leq n/5$ if and only if $G$ is not an $\mathcal{E}_4$-graph. They also conjectured \cite[Conjecture~3.4]{BBS} that for some finite set $\mathcal{E}$ of graphs, $\iota(G, \mathcal{C}') \leq n/5$ if and only if $G$ is not an $\mathcal{E}$-graph. Recently, Chen and Cui \cite{CC} verified this conjecture and proved that $\mathcal{E} = \mathcal{E}_4$. This had been proved by Zhang and Wu \cite{ZW4} for the case where $G$ contains no triangle. Suppose that $G$ is not a $4$-cycle. Wei, Zhang and Zhao \cite{WZZ} showed that 
\begin{equation} \iota(G, \mathcal{F}) \leq \frac{m+1}{6} \label{ZWineq}
\end{equation}
if $\mathcal{F} = \{C_4\}$. Zhang and Wu \cite{ZW3} showed that (\ref{ZWineq}) holds if $\mathcal{F} = \mathcal{C}$ and $G$ contains no triangle. For this result, we have $\iota(G,\mathcal{F}) = \iota(G,\mathcal{C}')$ due to the condition that $G$ contains no triangle. Generalizing both the Wei--Zhang--Zhao result and the Zhang--Wu result, we show that (\ref{ZWineq}) holds also if this condition is dropped and $\mathcal{F} = \mathcal{C}'$, hence providing an analogue of the Chen--Cui result. We also determine the extremal graphs. Let $C_4'$ be the \emph{diamond graph} $([4], E(C_4) \cup \{\{1,3\}\})$. We can now state our result, which is proved in the next section.

\begin{theorem} \label{mainresult}
If $G$ is a connected $m$-edge graph that is not a $4$-cycle, then
\begin{equation} \iota(G, \mathcal{C}') \leq \frac{m+1}{6}. \label{mainbd}
\end{equation} 
Moreover, the following statements hold:\\
(i) Equality in (\ref{mainbd}) holds if and only if $G$ is a pure $(m,C_4)$-special graph or a $\{C_4', C_5\}$-graph.\\
(ii) If $G$ is an $(m,C_4)$-special graph, then $\iota(G, \mathcal{C}') = \lfloor(m+1)/6 \rfloor$.
\end{theorem}

\noindent It is worth pointing out that the proof of (\ref{mainbd}) makes use of (i) in an inductive argument. In the literature on isolation, this idea of establishing a bound by also determining the graphs attaining it was introduced in \cite{Borg2}. In order to prove (\ref{ZWineq}), one can instead take the approach in \cite{BBS}, where the case where $G$ is subcubic is treated separately from the case where $G$ is not subcubic. This is the approach taken in \cite{WZZ}.

\section{Proof of Theorem~\ref{mainresult}} \label{Proofsection}

We start the proof of Theorem~\ref{mainresult} with two basic lemmas. 

\begin{lemma}[\cite{Borg1}] \label{lemma}
If $G$ is a graph, $\mathcal{F}$ is a set of graphs, $X \subseteq V(G)$ and $Y \subseteq N[X]$, then \[\iota(G, \mathcal{F}) \leq |X| + \iota(G-Y, \mathcal{F}).\] 
\end{lemma}

\begin{lemma}[\cite{Borg1, Borg4}] \label{lemmacomp} If $G_1, \dots, G_r$ are the distinct components of a graph $G$, and $\mathcal{F}$ is a set of connected graphs, then $\iota(G,\mathcal{F}) = \sum_{i=1}^r \iota(G_i,\mathcal{F})$.
\end{lemma}

The next lemma concerns a case where no member of a subset $Y$ of $V(G)$ is a vertex of an $\mathcal{F}$-graph contained by $G$, where $\mathcal{F}$ is a set of cycles.

\begin{lemma} [\cite{BBS}] \label{lemmaisolsquares}
If $G$ is a graph, $\mathcal{F}$ is a set of cycles, $x \in V(G)$, $Y \subseteq V(G) \backslash \{x\}$, $N[Y] \cap V(G-Y) \subseteq \{x\}$, and $G[\{x\} \cup Y]$ contains no $\mathcal{F}$-graph, then $\iota(G, \mathcal{F}) = \iota(G-Y, \mathcal{F})$ and every $\mathcal{F}$-isolating set of $G-Y$ is an $\mathcal{F}$-isolating set of $G$.
\end{lemma}

An \emph{isolated vertex of $G$} is a vertex of $G$ of degree 0. A \emph{leaf of $G$} is a vertex of $G$ of degree $1$.

\begin{cor} [\cite{BBS}] \label{Corollary}
If $G$ is a graph, $\mathcal{F}$ is a set of cycles, and $y$ is an isolated vertex of $G$ or a leaf of $G$, then $\iota(G, \mathcal{F})=\iota(G-y, \mathcal{F})$.
\end{cor}
Corollary~\ref{Corollary} generalizes as follows.

\begin{cor} \label{Corollary2}
If $G$ is a graph, $\mathcal{F}$ is a set of cycles, and $\emptyset \neq Y \subseteq V(G)$ such that each member of $Y$ is an isolated vertex of $G$ or a leaf of $G$, then $\iota(G, \mathcal{F}) = \iota(G-Y, \mathcal{F})$.
\end{cor}
\textbf{Proof.} We use induction on $|Y|$. If $|Y| = 1$, then the result is Corollary~\ref{Corollary}. Suppose $|Y| \geq 2$. Let $y \in Y$. By Corollary~\ref{Corollary}, $\iota(G-y, \mathcal{F}) = \iota(G, \mathcal{F})$. Let $Y' = Y \backslash \{y\}$. Then, each member of $Y'$ is an isolated vertex of $G-y$ or a leaf of $G-y$. By the induction hypothesis, $\iota(G-y, \mathcal{F}) = \iota((G-y)-Y', \mathcal{F})$. Since $\iota(G-y, \mathcal{F}) = \iota(G, \mathcal{F})$ and $\iota((G-y)-Y', \mathcal{F}) = \iota(G-Y, \mathcal{F})$, the result follows.~\hfill{$\Box$}
\\

For a vertex $v$ of a graph $G$, let $E_G(v)$ denote the set $\{vw \colon w \in N_G(v)\}$. For $X, Y \subseteq V(G)$, let $E_G(X,Y)$ denote the set $\{xy \in E(G) \colon x \in X, \, y \in Y\}$. Let ${\rm C}(G)$ denote the set of components of $G$.

\begin{lemma} \label{lemmaextension} If $G$ is a graph, $\mathcal{F}$ is a set of cycles, $Y \subseteq V(G)$, $G[Y]$ contains no $\mathcal{F}$-graph, and $|E_G(V(H), Y)| \leq 1$ for each $H \in {\rm C}(G-Y)$, then $\iota(G, \mathcal{F}) \leq \iota(G-Y, \mathcal{F})$ and every $\mathcal{F}$-isolating set of $G-Y$ is an $\mathcal{F}$-isolating set of $G$.
\end{lemma}

\noindent \textbf{Proof.} Suppose that $G$ contains an $\mathcal{F}$-graph $F$. Since $V(F) \nsubseteq V(G[Y])$, $V(F) \cap V(H) \neq \emptyset$ for some $H \in {\rm C}(G-Y)$. Let $X = V(H)$ and $Z = V(G) \backslash X$. Suppose $V(F) \nsubseteq X$. Then, $z_1z_2 \in E(F)$ for some $z_1 \in X$ and $z_2 \in Z$. Since $H \in {\rm C}(G-Y)$, $E_G(X, Z) = E_G(X, Y)$. Since $|E_G(X, Y)| \leq 1$, $z_2 = y$ for some $y \in Y$, and $E_G(X, Z) = \{z_1z_2\}$. We have $E(F) = \{z_1z_2, z_2z_3, \dots, z_{r-1}z_r, z_rz_1\}$ for some $r \geq 3$ and some distinct members $z_3, \dots, z_r$ of $V(G) \backslash \{z_1, z_2\}$. We have $N_F(y) = N_F(z_2) = \{z_1, z_3\}$. Since $E_G(X, Z) = \{z_1z_2\}$, $z_3 \notin X \cup \{y\}$. Since $y \notin \{z_3, \dots, z_r\}$ and $E_G(X, Z) = \{z_1z_2\}$, we obtain $z_3, \dots, z_r, z_1 \notin X \cup \{y\}$, which contradicts $z_1 \in X$. Therefore, $V(F) \subseteq X$. Consequently, if $D$ is an $\mathcal{F}$-isolating set of $G-Y$, then $\emptyset \neq N_{G-Y}[D] \cap V(F) \subseteq N_{G}[D] \cap V(F)$, meaning that $D$ is an $\mathcal{F}$-isolating set of $G$.~\hfill{$\Box$}
 
\begin{prop} \label{Proposition}
(a) If $G$ is a pure $(m,C_4)$-special graph with exactly $q$ $C_4$-constituents, then $m = 6q-1$, $\iota(G,\mathcal{C}') = q$, and for any $v \in V(G)$, $G$ has a $\mathcal{C}'$-isolating set $D$ with $v \in D$ and $|D| = q$. \\
(b) If $G$ is a $\{C_4', C_5\}$-graph, then $\iota(G,\mathcal{C}') = 1 = (|E(G)| + 1)/6$, and for any $v \in V(G)$, $\{v\}$ is a $\mathcal{C}'$-isolating set of $G$.
\end{prop}
\textbf{Proof.} Suppose that $G$ is a pure $(m,F)$-special graph with exactly $q$ $F$-constituents as in Construction~\ref{const}, where $F = C_4$. For some $j \in [q]$, $v \in V(G_j)$. Let $D = \{v\} \cup \{v_i \colon i \in [q] \backslash \{j\}\}$. Then, $D$ is a $\mathcal{C}'$-isolating set of $G$, so $\iota(G, \mathcal{C}') \leq q$. If $S$ is a $\mathcal{C}'$-isolating set of $G$, then, since $G_1 - v_1, \dots, G_q - v_q$ are $4$-cycles, $S \cap V(G_i) \neq \emptyset$ for each $i \in [q]$. Therefore, $\iota(G,\mathcal{C}') = q$. Now $m=5q + |E(T)|$. Since $T$ is a $q$-vertex tree, $|E(T)|=q-1$. Thus, $m=6q-1$. This settles (a). (b) is trivial.~\hfill{$\Box$}\\


\noindent
\textbf{Proof of Theorem~\ref{mainresult}.} The argument in the proof of Proposition~\ref{Proposition} yields (ii). Proposition~\ref{Proposition} settles the sufficiency condition in (i). We now prove (\ref{mainbd}) and (i), using induction on~$m$.


The result is trivial if $m \leq 4$ as $G \not\simeq C_4$. Suppose $m \geq 5$. Let $k$ be the maximum degree $\max\{d(v) \colon v \in V(G)\}$ of $G$. Since $G$ is connected, $k \geq 2$. If $k = 2$, then $G$ is a path or a cycle. If $G$ is a path, then $\iota(G, \mathcal{C}') = 0 < (m+1)/6$. If $G$ is a cycle, then $\iota(G,\mathcal{C}') = 1 \leq (m+1)/6$, and equality holds only if $m = 5$. Suppose $k \geq 3$. Let $v \in V(G)$ with $d(v) = k$. Suppose $N[v] = V(G)$. Then, $\{v\}$ is a $\mathcal{C}'$-isolating set of $G$, so $\iota(G, \mathcal{C}') \leq 1 \leq (m+1)/6$. If $\iota(G, \mathcal{C}') = (m+1)/6$, then $G \simeq C_4'$. Now suppose $V(G) \neq N[v]$. Let $G' = G-N[v]$. Then, $V(G') \neq \emptyset$. Let $\mathcal{H} = {\rm C}(G')$. Let $\mathcal{H}' = \{H \in \mathcal{H} \colon H \simeq C_4\}$. For each $H \in \mathcal{H} \backslash \mathcal{H}'$, let $D_H$ be a $\mathcal{C}'$-isolating set of $H$ of size $\iota(H,\mathcal{C}')$. By the induction hypothesis, for each $H \in \mathcal{H} \backslash \mathcal{H}'$, 
\[|D_H| \leq \frac{|E(H)|+1}{6},\] 
and equality holds only if $H$ is a pure $(|E(H)|,C_4)$-special graph or a $\{C_4', C_5\}$-graph. 

For any $H \in \mathcal{H}$ and any $x \in N(v)$ such that $xy_{x,H} \in E(G)$ for some $y_{x,H} \in V(H)$, we say that $H$ is \emph{linked to $x$} and that $x$ is \emph{linked to $H$}. Since $G$ is connected, each member $H$ of $\mathcal{H}$ is linked to at least one member of $N(v)$, so $x_Hy_H \in E(G)$ for some $x_H \in N(v)$ and some $y_H \in V(H)$. For each $x \in N(v)$, let $\mathcal{H}'_x = \{H \in \mathcal{H}' \colon H \mbox{ is linked to } x\}$ and $\mathcal{H}_x^* = \{H \in \mathcal{H} \backslash \mathcal{H}' \colon H \mbox{ is linked to $x$ only}\}$. 

Let $H_1, \dots, H_p$ be the distinct members of $\mathcal{H}$. For each $i \in [p]$, let $m_i = |E(H_i)|$. Let $I = \{i \in [p] \colon H_i \mbox{ is a pure $(m_i,C_4)$-special graph}\}$. For each $i \in I$, let $G_{i,1}, \dots, G_{i,q_i}$ be the distinct $C_4$-constituents of $H_i$, and for each $j \in [q_i]$, let $v_{i,j}$ be the $C_4$-connection of $G_{i,j}$ in $H_i$, and let $w^1_{i,j}, w_{i,j}^2, w_{i,j}^3, w^4_{i,j}$ be the members of $V(G_{i,j}) \backslash \{v_{i,j}\}$ such that $E(G_{i,j}) = \{v_{i,j}w^1_{i,j}, w_{i,j}^1 w^2_{i,j}, w_{i,j}^2 w^3_{i,j}, w_{i,j}^3 w^4_{i,j}, w_{i,j}^4 w^1_{i,j}\}$. By Proposition~\ref{Proposition}, $m_i + 1 = 6q_i$ for each $i \in I$. Let $I' = \{i \in [p] \colon H_i \mbox{ is a $\{C_4', C_5\}$-graph}\}$. Then, $m_i + 1 = 6$ for each $i \in I'$. Let $J = \{i \in [p] \colon H_i \simeq C_4\}$. Let $J' = [p] \backslash (I \cup I' \cup J)$. Then, $I$, $I'$, $J$ and $J'$ are pairwise disjoint, $\mathcal{H}' = \{H_i \colon i \in J\}$ and $\mathcal{H} \backslash \mathcal{H}' = \{H_i \colon i \in I \cup I' \cup J'\}$. If $i \in J'$, then by the induction hypothesis, $\iota(H_i, \mathcal{C}') < (m_i + 1)/6$, so $6\iota(H_i, \mathcal{C}') < m_i + 1$, and hence $6\iota(H_i, \mathcal{C}') \leq m_i$. Thus, 
\begin{equation} J' = \{i \in [p] \colon \iota(H_i, \mathcal{C}') \leq m_i/6\}. \label{J'eq}
\end{equation}
We have $|E_G(v)| = d(v) = k \geq 3$. Let $A_1 \subseteq E_G(v)$ with $|A_1| = 3$. Thus, $A_1 = \{vx_1, vx_2, vx_3\}$ for some $x_1, x_2, x_3 \in N(v)$. Let $A_2 = \{x_Hy_H \colon H \in \mathcal{H}\}$. Let $M_1 = E(G[N[v]])$, $M_2 = E_G(N(v), V(G'))$ and $M_3 = \bigcup_{H \in \mathcal{H}} E(H)$. Thus, $A_1 \subseteq M_1$, $A_2 \subseteq M_2$ and $m = |M_1| +|M_2| + |M_3|$. Let $a = |M_1 \backslash A_1| + |M_2 \backslash A_2|$ and $b = |\{x_{H_i}y_{H_i} \colon i \in J'\}|$. We have
%
%
\begin{align}
m &= |A_1| + |A_2| + a  + \sum_{i \in I \cup I' \cup J \cup J'} m_i \nonumber \\
&= 3 + a + b + \sum_{i \in J'} m_i + \sum_{i \in I \cup I' \cup J} |E(H_i) \cup \{x_{H_i}y_{H_i}\}|. \label{numberofedges} 
\end{align}

We need to treat the case $\mathcal{H}' = \emptyset$ differently from the case $\mathcal{H}' \neq \emptyset$, particularly because we cannot apply the induction hypothesis to members of $\mathcal{H}'$.
\\

\noindent \textbf{Case 1:} \emph{$\mathcal{H}' = \emptyset$.} Then, $J = \emptyset$. By Lemma~\ref{lemma} (with $X = \{v\}$ and $Y = N[v]$), Lemma~\ref{lemmacomp} and (\ref{numberofedges}),
\begin{align}
\iota(G, \mathcal{C}') &\leq 1 + \iota(G', \mathcal{C}') = 1 + \sum_{H \in \mathcal{H}} \iota(H, \mathcal{C}') \leq 1 + \sum_{i \in J'} \frac{m_i}{6} + \sum_{i \in I \cup I'} \frac{m_i + 1}{6} \nonumber \\
&= 1 + \frac{m - 3 - a - b}{6}= \frac{m + 3 - a - b}{6}. \nonumber
\end{align}
If $a+b \geq 3$, then $\iota(G, \mathcal{C}') \leq m/6$. Suppose $a+b \leq 2$. \\

\noindent \textbf{Case 1.1:} \emph{$G[N[v]]$ contains a $\mathcal{C}'$-graph $G_0$.} We have $|V(G_0)| \geq 4$. If we assume that $|N[v]| \geq 5$, then we obtain $a \geq 3$, which contradicts $a+b \leq 2$. Thus, $N[v] = \{v, x_1, x_2, x_3\} = V(G_0)$. We may assume that $E(G_0) = \{vx_1, x_1x_2, x_2x_3, x_3v\}$. Since $a+b \leq 2$, we have $E(G) = A_1 \cup \{x_1x_2, x_2x_3\} \cup A_2 \cup M_3$ and $b = 0$. Since $b = 0$, $J' = \emptyset$. By Proposition~\ref{Proposition}, for any $i \in I \cup I'$, $H_i$ has a $\mathcal{C}'$-isolating set $D_i$ with $y_{H_i} \in D_i$ and $|D_i| = (m_i + 1)/6$. 
Clearly, $\bigcup_{i \in I \cup I'} D_i$ is a $\mathcal{C}'$-isolating set of $G$, so $\iota(G, \mathcal{C}') \leq \sum_{i \in I \cup I'} |E(H_i) \cup \{x_{H_i}y_{H_i}\}|/6 = (m-5)/6$.
\\

\noindent \textbf{Case 1.2:} \emph{$G[N[v]]$ contains no $\mathcal{C}'$-graph.} Suppose $M_2 = A_2$. Then, $E_G(V(H), N[v]) = \{x_Hy_H\}$ for each $H \in \mathcal{H}$. By Lemmas~\ref{lemmaextension} and~\ref{lemmacomp},
\[\iota(G, \mathcal{C}') \leq \iota(G', \mathcal{C}') = \sum_{H \in \mathcal{H}} \iota(H, \mathcal{C}') \leq \sum_{H \in \mathcal{H}} \frac{|E(H) \cup \{x_Hy_H\}|}{6} = \frac{m-3-a}{6}.\bigskip\]
Now suppose $M_2 \neq A_2$. Then, $a \geq 1$ and there exist $x \in N(v)$, $i \in [p]$ and $y \in V(H_i)$ such that $xy \in E(G)$ and $xy \neq x_{H_i}y_{H_i}$. Let $G^* = G - V(H_i)$. Clearly, $G^*$ is connected and is not a $4$-cycle (as $|N[v]| \geq 4$ and $G[N[v]] \not\simeq C_4$). By the induction hypothesis,
\[\iota(G^*, \mathcal{C}') \leq \frac{(m - |E(H_i) \cup \{x_{H_i}y_{H_i}, xy\}|) + 1}{6} = \frac{m - m_i - 1}{6}.\]
Let $D^*$ be a $\mathcal{C}'$-isolating set of $G^*$ of size $\iota(G^*, \mathcal{C}')$.
\\

\noindent \textbf{Case 1.2.1:} \emph{$i \in I$.} Let $X = \{w_{i,j}^1 \colon j \in [q_i]\}$, $Y = \{w_{i,j}^3 \colon j \in [q_i]\}$ and $Y' = V(H_i) \backslash Y$. Let $G_Y = G-Y'$. Suppose $d_{G_Y}(w) \leq 1$ for each $w \in Y$. Then, by Corollary~\ref{Corollary2}, $\iota(G_Y, \mathcal{C}') = \iota(G_Y - Y, \mathcal{C}') = \iota(G^*, \mathcal{C}')$. Since $Y' \subseteq N[X]$, Lemma~\ref{lemma} gives us 
\[\iota(G, \mathcal{C}') \leq |X| + \iota(G_Y, \mathcal{C}') = q_i + \iota(G^*, \mathcal{C}') \leq \frac{m_i + 1}{6} + \frac{m - m_i - 1}{6} = \frac{m}{6}.\]
Now suppose $d_{G_Y}(w) \geq 2$ for some $w \in Y$. We have $w = w_{i,j}^3$ for some $j \in [q_i]$. Since $w_{i,j}^2, w_{i,j}^4 \in N(w) \backslash V(G_Y)$, $d(w) \geq 4$. Since $d(v) = \Delta(G)$, $d(v) \geq 4$. Since $a+b \leq 2$, it follows that $M_1 \backslash A_1 = \{vx_4\}$ for some $x_4 \notin \{x_1, x_2, x_3\}$, $M_2 \backslash A_2 = \{xy\}$, $E_{G_Y}(w) = \{x_{H_i}y_{H_i}, xy\} = \{x_{H_i} w, xw\}$ and $E(G) = A_1 \cup A_2 \cup \{vx_4, xw\} \cup M_3$. Let $X' = \{w\} \cup (X \backslash \{w_{i,j}^1\})$ and $D = X' \cup D^*$. Since $E_G(V(G^*), V(H_i)) \subseteq E_G(w)$ and $X'$ is a $\mathcal{C}'$-isolating set of $H_i$, $D$ is a $\mathcal{C}'$-isolating set of $G$. Therefore, we have
\[\iota(G, \mathcal{C}') \leq |X'| + |D^*| = q_i + \iota(G^*, \mathcal{C}') \leq \frac{m_i + 1}{6} + \frac{m - m_i - 1}{6} = \frac{m}{6}.\]
\medskip

\noindent \textbf{Case 1.2.2:} \emph{$i \in I'$.} Then, $H_i \simeq C_4'$ or $H_i \simeq C_5$. Thus, there exists some $w \in V(H_i)$ such that $y_{H_i}, y \in N[w]$ (because if $y_{H_i} \neq y$, then $y_{H_i}$ and $y$ are of distance at most $2$ in $H_i$). 
Let $Y = V(H_i) \backslash N_{H_i}[w]$, $Y' = N_{H_i}[w]$ and $G_Y = G - Y'$. Then, $x_{H_i}y_{H_i}, xy \notin E(G_Y)$. Since $a + b \leq 2$, $|E_{G_Y}(V(G^*), Y)| \leq 1$. Thus, for some $x^* \in V(G^*)$, $N_{G_Y}[Y] \cap V(G^*) \subseteq \{x^*\}$ and $G_Y[\{x^*\} \cup Y]$ contains no $\mathcal{C}'$-graph. Since $G^* = G_Y - Y$, Lemma~\ref{lemmaisolsquares} yields $\iota(G_Y, \mathcal{C}') = \iota(G^*, \mathcal{C}')$. By Lemma~\ref{lemma},
\[\iota(G, \mathcal{C}') \leq 1 + \iota(G_Y, \mathcal{C}') = 1 + \iota(G^*, \mathcal{C}') \leq \frac{m_i + 1}{6} + \frac{m - m_i - 1}{6} = \frac{m}{6}.\]
\medskip

\noindent \textbf{Case 1.2.3:} \emph{$i \in J'$.} Then, $b \geq 1$. Since $a+b \leq 2$ and $a \geq 1$, we have $a = b = 1$, $J' = \{i\}$ and $E(G) = A_1 \cup A_2 \cup \{xy\} \cup M_3$. For each $j \in I'$, let $D_j = \{y_{H_j}\}$. By Proposition~\ref{Proposition}, for each $j \in I$, $H_j$ has a $\mathcal{C}'$-isolating set $D_j$ with $y_{H_j} \in D_j$ and $|D_j| = q_j$. Let $X = \bigcup_{j \in [p] \backslash \{i\}} (N[D_j] \cap V(H_j))$ and $D_X = \bigcup_{j \in [p] \backslash \{i\}} D_j$. Let $G_v^* = G[N[v] \cup V(H_i)]$. Then, $G_v^*$ is a component of $G-X$, and any other component of $G-X$ contains no $\mathcal{C}'$-graph. Let $D_v^*$ be a $\mathcal{C}'$-isolating set of $G_v^*$ of size $\iota(G_v^*, \mathcal{C}')$. Then, $D_v^* \cup D_X$ is a $\mathcal{C}'$-isolating set of $G$. Let $x' \in \{x_1, x_2, x_3\} \backslash \{x_{H_i}, x\}$. Since $E_G(v) = A_1$ and $E_G(N(v), V(H_i)) = \{x_{H_i}y_{H_i}, xy\}$, $x'$ is a leaf of $G_v^*$. By Corollary~\ref{Corollary}, $\iota(G_v^*, \mathcal{C}') = \iota(G_v^* - x', \mathcal{C}')$. 

Suppose $G_v^* - x' \not\simeq C_4$. By the induction hypothesis, $\iota(G_v^* - x', \mathcal{C}') \leq (|E(G_v^* - x')| + 1)/6$, so $\iota(G_v^*, \mathcal{C}') \leq |E(G_v^*)|/6$. We have
\begin{align} \iota(G, \mathcal{C}') &\leq |D_v^*| + |D_X| = \iota(G_v^*, \mathcal{C}') + \sum_{j \in [p] \backslash \{i\}} |D_j| \nonumber \\
&\leq \frac{|E(G_v^*)|}{6} + \sum_{j \in [p] \backslash \{i\}} \frac{|E(H_j) \cup \{x_{H_j}y_{H_j}\}|}{6} = \frac{m}{6}. \nonumber
\end{align}

Now suppose $G_v^* - x' \simeq C_4$. Then, $V(H_i) = \{y\} = \{y_{H_i}\}$, $x \neq x_{H_i}$ and $E(G_v^* - x') = \{vx_{H_i}, x_{H_i}y, yx, xv\}$. If $I \cup I' = \emptyset$, then $E(G) = E(G_v^* - x') \cup \{vx'\}$, so $G$ is a pure $(5, C_4)$-special graph. Suppose $I \cup I' \neq \emptyset$.

Suppose $I' \neq \emptyset$. Let $h \in I'$. Then, $(D_X \backslash D_h) \cup \{x_{H_h}\}$ is a $\mathcal{C}'$-isolating set of $G$, so 
\begin{equation} \iota(G, \mathcal{C}') \leq |D_X| = \sum_{j \in [p] \backslash \{i\}} |D_j| = \sum_{j \in [p] \backslash \{i\}} \frac{|E(H_j) \cup \{x_{H_j}y_{H_j}\}|}{6} < \frac{m}{6}. \nonumber
\end{equation}

Now suppose $I' = \emptyset$. Then, $I \neq \emptyset$. Suppose $y_{H_h} \notin \{v_{h,j} \colon j \in [q_h]\}$ for some $h \in I$. Then, $y_{H_h} = w_{h,j'}^t$ for some $j' \in [q_h]$ and $t \in [4]$. Let $D_h' = \{x_{H_h}\} \cup \{v_{h,j} \colon j \in [q_h] \backslash \{j'\}\}$. Then, $(D_X \backslash D_h) \cup D_h'$ is a $\mathcal{C}'$-isolating set of $G$, so $\iota(G, \mathcal{C}') < m/6$ as above. Now suppose $y_{H_h} \in \{v_{h,j} \colon j \in [q_h]\}$ for each $h \in I$. If $x_{H_h} \neq x'$ for some $h \in I$, then $D_X$ is a $\mathcal{C}'$-isolating set of $G$, so $\iota(G, \mathcal{C}') < m/6$ as above. Suppose $x_{H_h} = x'$ for each $h \in I$. Then, $G[\{v_{h,j} \colon h \in I, \, j \in [q_h]\} \cup \{x'\}]$ is a tree, and $G$ is a pure $(m,C_4)$-special graph whose $C_4$-constituents are $G[N[v] \cup \{y\}]$ and the members of $\{G_{h,j} \colon h \in I, \, j \in [q_h]\}$. 
\\

\noindent\textbf{Case 2:} \emph{$\mathcal{H}'\neq \emptyset$.} Let $H' \in \mathcal{H}'$. Let $x \in N(v)$ such that $H'$ is linked to $x$. Thus, $H'$ is a $4$-cycle $(\{y_1,y_2,y_3,y_4\}, \{y_1y_2, y_2y_3, y_3y_4, y_4y_1\})$ with $xy_1 \in E(G)$.  Let $\mathcal{H}_1 = \{H \in \mathcal{H}' \colon H \mbox{ is linked to $x$ only}\}$ and $\mathcal{H}_2 = \{H \in \mathcal{H}\backslash \mathcal{H}' \colon H \mbox{ is linked to $x$ only}\}$.
\\ 

\noindent \textbf{Case 2.1:} \emph{Each member of $\mathcal{H}'$ is linked to at least two members of $N(v)$.} Then, $\mathcal{H}_1 = \emptyset$ and $H'$ is linked to some $x' \in N(v) \backslash \{x\}$, so $x'y' \in E(G)$ for some $y'\in V(H')$. Let $Y = \{x, y_1, y_2, y_4\}$ and $G^* = G-Y$. Then, $G^*$ has a component $G_v^*$ with $N[v]\backslash\{x\} \subseteq V(G_v^*)$, and $\{G_v^*\} \cup \mathcal{H}_2 \subseteq {\rm C}(G^*) \subseteq \{G_v^*, (\{y_3\}, \emptyset)\} \cup \mathcal{H}_2$. 
If $y_3 \notin V(G_v^*)$, then $y_3$ is an isolated vertex of $G^*$, so $\iota(G^*, \mathcal{C}') = \iota(G^*-y_3, \mathcal{C}')$ by Corollary~\ref{Corollary}. By Lemma~\ref{lemmacomp}, $\iota(G^*, \mathcal{C}') = \iota(G_v^*, \mathcal{C}') + \sum_{H \in \mathcal{H}_2} \iota(H, \mathcal{C}')$. Since $Y \subseteq N[y_1]$, Lemma~\ref{lemma} yields
\[\iota(G, \mathcal{C}') \leq 1 + \iota(G^*, \mathcal{C}') \leq \frac{|E(H') \cup \{xy_1, vx\}|}{6} + \iota(G_v^*, \mathcal{C}') + \sum_{H \in \mathcal{H}_2} \frac{|E(H) \cup \{xy_{x,H}\}|}{6}.\]

Suppose $G_v^* \not\simeq C_4$. By the induction hypothesis, $\iota(G_v^*, \mathcal{C}') \leq (|E(G_v^*)| + 1)/6$, so
$\iota(G, \mathcal{C}') \leq (m+1)/6$. Suppose $\iota(G, \mathcal{C}')= (m+1)/6$. Then, $\iota(G, \mathcal{C}') = 1 + \iota(G^*, \mathcal{C}')$, $\iota(G_v^*, \mathcal{C}')=(|E(G_v^*)|+1)/6$, $\iota(H, \mathcal{C}')=(|E(H)|+1)/6$ for each $H \in \mathcal{H}_2$, and
\begin{equation} E(G) = E(H') \cup \{xy_1, vx\} \cup E(G_v^*) \cup \bigcup_{H \in \mathcal{H}_2} (E(H) \cup \{xy_{x,H}\}). \label{2.1eqn}
\end{equation}
%
By the induction hypothesis, each member $F$ of $\{G_v^*\} \cup \mathcal{H}_2$ is a pure $(|E(F)|, C_4)$-special graph or a $\{C_4', C_5\}$-graph. By Proposition~\ref{Proposition}, $G_v^*$ has a $\mathcal{C}'$-isolating set $D_v^*$ with $x' \in D_v^*$ and $|D_v^*| = \iota(G_v^*, \mathcal{C}')$, and for each $H \in \mathcal{H}_2$, $H$ has a $\mathcal{C}'$-isolating set $D_H'$ with $y_{x,H} \in D_H'$ and $|D_H'| = \iota(H, \mathcal{C}')$. Let $D = D_v^* \cup \bigcup_{H \in \mathcal{H}_2} D_H'$. By (\ref{2.1eqn}), we have $x'y' \in E(G_v^*)$, so $y' = y_3 \in V(G_v^*)$. 
Also by (\ref{2.1eqn}), $E_G(V(G_v^*), Y) = \{vx, y_2y_3, y_3y_4\}$ and $E(G[Y]) = \{xy_1, y_1y_2, y_1y_4\}$. Let $H'' = G[Y]$ if $v \notin D$ and $\mathcal{H}_2 = \emptyset$, and let $H'' = G[Y] - x$ if $v \in D$ or $\mathcal{H}_2 \neq \emptyset$. Note that $x \in N[D]$ if and only if $v \in D$ or $\mathcal{H}_2 \neq \emptyset$. Since $v, y_3 \in N(x')$, ${\rm C}(G-N[D]) = \{H''\} \cup {\rm C}(G_v^* - N_{G_v^*}[D_v^*]) \cup \bigcup_{H \in \mathcal{H}_2} {\rm C}(H - N_H[D_H'])$. Thus, we have that $D$ is a $\mathcal{C}'$-isolating set of $G$ of size $\iota(G^*, \mathcal{C}')$, contradicting $\iota(G, \mathcal{C}') = 1 + \iota(G^*, \mathcal{C}')$. 

Now suppose $G_v^* \simeq C_4$. Then, $\mathcal{H}' = \{H'\}$ and $d(v) = 3 = \Delta(G)$. We may assume that $x = x_1$ and $x' = x_2$. We have $G_v^* = (\{v, x_2, x_3, z\}, \{vx_2, x_2z, zx_3, x_3v\})$ for some $z \in V(G') \backslash Y$. Since $\Delta(G) = 3$, we obtain $z \notin V(H')$ (otherwise, $d(z) \geq |\{x_2, x_3\} \cup N_{H'}(z)| = 4$), $N(y_3) \backslash \{y_2, y_4\} \subseteq \{x_1\}$ (because if $N(y_3) \cap \{x_2, x_3\} \neq \emptyset$, then $y_3 \in V(G_v^*)$), $N(y_1) = \{x_1, y_2, y_4\}$ and $y' \in \{y_2,y_4\}$. We may assume that $y' = y_2$. Thus, $N(x_2) = \{v, z, y_2\}$, $N(y_2) = \{y_1, y_3, x_2\}$, $N(x_3) \backslash \{v, z\} \subseteq \{x_1, y_4\}$ and $N(y_4) \backslash \{y_1, y_3\} \subseteq \{x_1, x_3\}$. Let $X_1 = \{v, x_1, y_1\}$ and $F_1 = G[\{x_3, z, x_2, y_2, y_3, y_4\}]$. Then, ${\rm C}(G - X_1) = \{F_1\} \cup \mathcal{H}_2$. Since $X_1 \subseteq N[x_1]$, $\iota(G, \mathcal{C}') \leq 1 + \iota(G-X_1, \mathcal{C}')$ by Lemma~\ref{lemma}. If $x_3y_4 \notin E(G)$, then $F_1$ is a path, so $\bigcup_{H \in \mathcal{H}_2} D_H$ is a $\mathcal{C}'$-isolating set of $G-X_1$, and hence we have
\[\iota(G, \mathcal{C}') \leq 1 + \sum_{H \in \mathcal{H}_2} \frac{|E(H)|+1}{6} < \frac{|E(H')\cup E(G_v^*)|}{6} + \sum_{H \in \mathcal{H}_2} \frac{|E(H) \cup \{xy_{x,H}\}|}{6} < \frac{m}{6}.\]
Suppose $x_3y_4 \in E(G)$. Let $X_2 = \{v, x_1, x_3, z, y_1, y_2, y_4\}$. Then, $\bigcup_{H \in \mathcal{H}_2} D_H$ is a $\mathcal{C}'$-isolating set of $G - X_2$ as ${\rm C}(G - X_2) = \{(\{x_2\}, \emptyset), (\{y_3\}, \emptyset)\} \cup \mathcal{H}_2$. Since $X_2 = N[\{y_1, x_3\}]$, Lemma~\ref{lemma} yields $\iota(G, \mathcal{C}') \leq 2 + \iota(G-X_2, \mathcal{C}')$, so
\begin{align} \iota(G, \mathcal{C}') &\leq \frac{|E(H') \cup E(G_v^*) \cup \{vx_1,x_1y_1,x_2y_2,x_3y_4\}|}{6} + \sum_{H \in \mathcal{H}_2} \frac{|E(H)  \cup \{xy_{x,H}\}|}{6} \leq \frac{m}{6}. \nonumber
\end{align}

\noindent \textbf{Case 2.2:} \emph{Some member of $\mathcal{H}'$ is linked to only one member of $N(v)$.} Thus, we may assume that $H'$ is linked to $x$ only. Let $h_1 = |\mathcal{H}_1|$. Then, $h_1 \geq 1$ as $H' \in \mathcal{H}_1$. 
Let $X = \{x\} \cup \bigcup_{H \in \mathcal{H}_1} V(H)$ and $D_X = \{x\}$. Then, $D_X$ is a $\mathcal{C}'$-isolating set of $G[X]$, and
\begin{equation}
|D_X| = 1 \leq \frac{1 + \sum_{H \in \mathcal{H}_1} |E(H) \cup \{xy_{x,H} \}|}{6} = \frac{5h_1 + 1}{6}. \nonumber
\end{equation}
Let $G^* = G-X$. Then, $G^*$ has a component $G^*_v$ with $N[v] \backslash \{x\} \subseteq V(G^*_v)$, and ${\rm C}(G^*) = \{G_v^*\} \cup \mathcal{H}_2$. 

If $G^*_v \not\simeq C_4$, then by the induction hypothesis, $G^*_v$ has a $\mathcal{C}'$-isolating set $D_v^*$ with $|D_v^*| = \iota(G_v^*, \mathcal{C}') \leq (|E(G^*_v)| + 1)/6$. If $G^*_v \simeq C_4$, then let $D_v^* = \{x\}$. Let $D = D^*_v \cup D_X \cup \bigcup_{H \in \mathcal{H}_2} D_H$. By the definition of $\mathcal{H}_1$ and $\mathcal{H}_2$, ${\rm C}(G-x) = \{G^*_v\} \cup \mathcal{H}_1 \cup \mathcal{H}_2$. Thus, $D$ is a $\mathcal{C}'$-isolating set of $G$ as $x \in D$, $v \in V(G^*_v) \cap N[x]$ and $D_X$ is a $\mathcal{C}'$-isolating set of $G[X]$. We have
\begin{align}
m &\geq |E(G^*_v) \cup \{vx\}| + \sum_{H \in \mathcal{H}_1 \cup \mathcal{H}_2} |E(H) \cup \{xy_{x,H}\}|  \nonumber \\
&= |E(G^*_v)| + 1 + 5h_1 + \sum_{H \in \mathcal{H}_2} (|E(H)|+1). \label{edgesineq2}
\end{align}

Suppose $G^*_v \simeq C_4$. Then, $D = \{x\} \cup \bigcup_{H \in \mathcal{H}_2} D_H$ and, by (\ref{edgesineq2}), 
\[m \geq 5(h_1 + 1) + \sum_{H \in \mathcal{H}_2} (|E(H)|+1) \geq 10 + \sum_{H \in \mathcal{H}_2} (|E(H)|+1)\] 
as $h_1 \geq 1$. We have
\begin{align}
\iota(G, \mathcal{C}') &\leq |D| = 1 + \sum_{H \in \mathcal{H}_2} |D_H| < \frac{10}{6} + \sum_{H \in \mathcal{H}_2} \frac{|E(H)| + 1}{6} \leq \frac{m}{6}. \nonumber
\end{align}

Now suppose $G^*_v \not\simeq C_4$. We have
\begin{align}\label{edgesboundineqsharp4}
\iota(G, \mathcal{C}') &\leq |D| = |D_v^*| + |D_X| + \sum_{H \in \mathcal{H}_2} |D_H| \leq \frac{|E(G^*_v)| + 1}{6} + 1 + \sum_{H \in \mathcal{H}_2} \frac{|E(H)| + 1}{6} \nonumber \\
&\leq \frac{|E(G^*_v)|+1}{6} +\frac{5h_1 + 1}{6} + \sum_{H \in \mathcal{H}_2} \frac{|E(H)|+1}{6} \leq \frac{m+1}{6} \quad \mbox{(by  (\ref{edgesineq2}))}.
\end{align}
Suppose $\iota(G, \mathcal{C}') = (m+1)/6$. Then, equality holds throughout in each of (\ref{edgesineq2}) and (\ref{edgesboundineqsharp4}). Consequently, $\iota(G, \mathcal{C}') = |D|$, $h_1 = 1$, $\mathcal{H}_1 = \{H'\}$, $|D_v^*| = (|E(G_v^*)| + 1)/6$, $|D_H| = (|E(H)| + 1)/6$ for each $H \in \mathcal{H}_2$, and
\[E(G) = E(G^*_v) \cup \{vx\} \cup E(H') \cup \{xy_{x,H'}\} \cup \bigcup_{H \in \mathcal{H}_2} (E(H) \cup \{xy_{x,H}\}).\]
Let $\mathcal{I} = \{F \in \{G_v^*\} \cup \mathcal{H}_2 \colon F \mbox{ is a pure $(|E(F)|, C_4)$-special graph}\}$ and $\mathcal{I}' = \{F \in \{G_v^*\} \cup \mathcal{H}_2 \colon F \mbox{ is a $\{C_4', C_5\}$-graph}\}$. By the induction hypothesis, $\{G_v^*\} \cup \mathcal{H}_2 = \mathcal{I} \cup \mathcal{I}'$. 

Suppose $\mathcal{I}' \neq \emptyset$. Let $F' \in \mathcal{I}'$. If $F' = G_v^*$, then let $D' = D \backslash D_v^*$. If $F' \in \mathcal{H}_2$, then let $D' = D \backslash D_{F'}$. Since $x \in D'$, $D'$ is a $\mathcal{C}'$-isolating set of $G$. We have $\iota(G, \mathcal{C}') \leq |D| - 1$, a contradiction. Therefore, $\mathcal{I}' = \emptyset$, and hence $\mathcal{I} = \{G_v^*\} \cup \mathcal{H}_2$. 

Let $F_1, \dots, F_r$ be the distinct members of $\mathcal{I}$, where $F_1 = G_v^*$. Let $y_{x, F_1} = v$ and $D_{F_1} = D_v^*$. For each $i \in [r]$, let $F_{i,1}, \dots, F_{i,s_i}$ be the distinct $C_4$-constituents of $F_i$, and for each $j \in [s_i]$, let $u_{i,j}^0$ be the $C_4$-connection of $F_{i,j}$, and let $u_{i,j}^1, u_{i,j}^2, u_{i,j}^3, u_{i,j}^4$ be the members of $V(F_{i,j}) \backslash \{u_{i,j}^0\}$ such that $E(F_{i,j}) = \{u_{i,j}^0u_{i,j}^1, u_{i,j}^1u_{i,j}^2, u_{i,j}^2u_{i,j}^3, u_{i,j}^3u_{i,j}^4, u_{i,j}^4u_{i,j}^1\}$. For each $i \in [r]$, let $D_i = \{u_{i,j}^0 \colon j \in [s_i]\}$. Then, $D_i$ is a $\mathcal{C}'$-isolating set of $F_i$. Since $|D_{F_i}| = (|E(F_i)| + 1)/6 = s_i$, $|D_{F_i}| = |D_i|$. Let $D' = \{x\} \cup \bigcup_{i=1}^r D_i$. Then, $D'$ is a $\mathcal{C}'$-isolating set of $G$ (as $x \in D'$) and $|D'| = |D|$. Suppose $y_{x, F_h} \notin D_h$ for some $h \in [r]$. Then, $y_{x, F_h} = u_{h,j'}^t$ for some $j' \in [s_h]$ and $t \in [4]$. Since $y_{x, F_h} \in N(x)$, $D' \backslash \{u_{h,j'}^0\}$ is a $\mathcal{C}'$-isolating set of $G$. We have $\iota(G, \mathcal{C}') \leq |D| - 1$, a contradiction. Thus, $y_{x, F_i} \in D_i$ for each $i \in [r]$. Therefore, $G[D']$ is a tree, and $G$ is a pure $(m,C_4)$-special graph whose $C_4$-constituents are $G[V(H') \cup \{x\}]$ and the members of $\{F_{i,j} \colon i \in [r], \, j \in [s_i]\}$.~\hfill{$\Box$}
\\
\\
\textbf{Acknowledgements:} The authors are grateful to the anonymous referee for checking the paper and providing constructive remarks.


\footnotesize
\begin{thebibliography}{}
\bibitem{BBS} K. Bartolo, P. Borg, D. Scicluna, Isolation of squares in graphs, Discrete Mathematics 347 (2024), paper 114161.

\bibitem{Borg1} P. Borg, Isolation of cycles, Graphs and Combinatorics 36 (2020), 631--637.

\bibitem{Borg2} P. Borg, Isolation of regular graphs and $k$-chromatic graphs, Mediterranean Journal of Mathematics 21 (2024), paper 148.

\bibitem{Borg3} P. Borg, Proof of a conjecture on isolation of graphs dominated by a vertex, Discrete Applied Mathematics 371 (2025), 247--253.

\bibitem{Borg4} P. Borg, Isolation of regular graphs, stars and $k$-chromatic graphs, Discrete Mathematics 349 (2026), paper 114706.

\bibitem{BFK} P. Borg, K. Fenech, P. Kaemawichanurat, Isolation of $k$-cliques, Discrete Mathematics 343 (2020), paper 111879.

\bibitem{BFK2} P. Borg, K. Fenech, P. Kaemawichanurat, Isolation of $k$-cliques II, Discrete Mathematics 345 (2022), paper 112641.

\bibitem{BCHH} R. Boutrig, M. Chellali, T.W. Haynes, S.T. Hedetniemi, Vertex-edge domination in graphs, Aequationes Mathematicae 90 (2016), 355--366.

\bibitem{BG} G. Boyer, W. Goddard, Disjoint isolating sets and graphs with maximum isolation number, Discrete Applied Mathematics 356 (2024), 110--116.

\bibitem{CaHa17} Y. Caro, A. Hansberg, Partial domination - the isolation number of a graph, Filomat 31 (2017), 3925--3944. 

\bibitem{CC} S. Chen, Q. Cui, Proof of a conjecture on isolation of cycles in graphs, Discrete Mathematics 349 (2026), 115076.

\bibitem{CCZ} S. Chen, Q. Cui, J. Zhang, A characterization of graphs with maximum cycle isolation number, Discrete Applied Mathematics 366 (2025), 161--175.

\bibitem{CCZ2} S. Chen, Q. Cui, L. Zhong, A characterization of graphs with maximum $k$-clique isolation number, Discrete Mathematics 348 (2025), 114531.

\bibitem{C} E.J. Cockayne, Domination of undirected graphs -- A survey, in: Lecture Notes in Mathematics, Volume 642, Springer, 1978, 141--147.

\bibitem{CH} E.J. Cockayne, S.T. Hedetniemi, Towards a theory of domination in graphs, Networks 7 (1977), 247--261.

\bibitem{CZ} Q. Cui, J. Zhang, A sharp upper bound on the cycle isolation number of graphs, Graphs and Combinatorics 39 (2023), paper 117.

\bibitem{HHS} T.W. Haynes, S.T. Hedetniemi, P.J. Slater, Fundamentals of Domination in Graphs, Marcel Dekker, Inc., New York, 1998.

\bibitem{HHS2} T.W. Haynes, S.T. Hedetniemi, P.J. Slater (Editors), Domination in Graphs: Advanced Topics, Marcel Dekker, Inc., New York, 1998.

\bibitem{HL2} S.T. Hedetniemi, R.C. Laskar, Bibliography on domination in graphs and some basic definitions of domination parameters, Discrete Mathematics 86 (1990), 257--277.

\bibitem{HL} S.T. Hedetniemi, R.C. Laskar (Editors), Topics on Domination, in: Annals of Discrete Mathematics, Volume 48, North-Holland Publishing Co., Amsterdam, 1991, Reprint of Discrete Mathematics 86 (1990).

\bibitem{LMS} M. Lema\'{n}ska,  M. Mora, M.J. Souto-Salorio, Graphs with isolation number equal to one third of the order, Discrete Mathematics 347 (2024), paper 113903.

\bibitem{Ore} O. Ore, Theory of graphs, American Mathematical Society Colloquium Publications, Volume 38, American Mathematical Society, Providence, R.I., 1962.

\bibitem{WZZ} X. Wei, G. Zhang, B. Zhao, On the $C_4$-isolation number of a graph, Discrete Mathematics, Algorithms and Applications, https://doi.org/10.1142/S1793830926500679. 

\bibitem{West} D.B. West, Introduction to Graph Theory, second edition, Prentice Hall, 2001.

\bibitem{Y} J. Yan, Isolation of the diamond graph, Bulletin of the Malaysian Mathematical Sciences Society 45 (2022), 1169--1181.

\bibitem{ZW3} G. Zhang, B. Wu, Cycle isolation of graphs with small girth, Graphs and Combinatorics 40 (2024), paper 38.

\bibitem{ZW2} G. Zhang, B. Wu, A note on the cycle isolation number of graphs,  Bulletin of the Malaysian Mathematical Sciences Society  47 (2024), paper 57.

\bibitem{ZW4} G. Zhang, B. Wu, On the cycle isolation number of triangle-free graphs, Discrete Mathematics 347 (2024), 114190.

\bibitem{Z} P. \.{Z}yli\'{n}ski, Vertex-edge domination in graphs, Aequationes Mathematicae 93 (2019), 735--742.
	
\end{thebibliography}
\end{document}